\documentclass{amsproc}

\input{epsf}

\newtheorem{theorem}{Theorem}[section]
\newtheorem{lemma}[theorem]{Lemma}
\newtheorem{proposition}[theorem]{Proposition}
\newtheorem{corollary}[theorem]{Corollary}

\theoremstyle{definition}
\newtheorem{definition}[theorem]{Definition}
\newtheorem{example}[theorem]{Example}

\theoremstyle{remark}
\newtheorem{remark}[theorem]{Remark}

\begin{document}

\title[Quandles and Cocycle Knot Invariants]{Diagrammatic Computations for 
Quandles and Cocycle Knot Invariants}

\author{J. Scott Carter}
\address{University of South Alabama,
Mobile, AL 36688}
\email{carter@mathstat.usouthal.edu}
\thanks{The first author was supported in part by NSF Grant \#9988107.}

\author{Seiichi Kamada}
\address{Osaka City University,
Osaka 558-8585, JAPAN}
\email{kamada@sci.osaka-cu.ac.jp }
\thanks{The second  author was supported by
  Fellowships from the Japan Society for the Promotion of Science.}

\author{Masahico Saito}
\address{ University of South Florida,
  Tampa, FL 33620 }
\email{saito@math.usf.edu}
\thanks{The third author was supported in part by NSF Grant \#9988101.}

\date{February 1st, 2001}

\begin{abstract}
The state-sum invariants for knots and knotted surfaces defined
from quandle cocycles are described  using the Kronecker product
between cycles  represented by colored knot diagrams and
a cocycle of a finite quandle used to color the diagram.
Such an interpretation is applied to evaluating the invariants.

Algebraic interpretations of quandle cocycles as deformations of
extensions are also given. The proofs rely on colored knot diagrams.
\end{abstract}

\maketitle

\section{Introduction}

Diagrammatic morphisms
interconnect algebra and topology.
Complicated algebraic formulas can be
  established via topological diagrams,
and algebraic structures give topological invariants.
In this paper, we present two
instances of algebraic and topological interplay
from
quandle homology theory.
First, we use the Kronecker product and
computations on colored knot diagrams for
evaluating the quandle knot cocycle
invariants. Second, we describe extensions of quandles by cocycles,
and give diagrammatics proofs. In both examples, relations
between algebra and diagrams play key roles.

A quandle is a set with a
self-distributive binary operation (defined below)
whose definition was motivated from knot theory.
A (co)homology theory was defined in \cite{CJKLS} for quandles,
which is a modification of rack (co)homology defined in \cite{FRS2}.
State-sum invariants using quandle cocycles as weights are
defined \cite{CJKLS} and computed for important families
of classical knots and knotted surfaces \cite{CJKS1}.
Quandle homomorphisms and virtual knots are applied to this
homology theory \cite{CJKS2}.
  The invariants were applied to study
knots,
for example, in detecting non-invertible
knotted surfaces \cite{CJKLS}.
On the other hand, knot diagrams colored by quandles can be used
to study quandle homology groups. This view point was developed
  in \cite{FRS2,Flower,Greene}
for rack homology and homotopy and generalized to quandle homology
in \cite{SSS2}.
It was poined out by Fenn and  Rourke that
the state-sum  terms can be interpreted as
Kronecker products. In this paper, we use such  interpretations
   to evaluate the invariants.

The second diagrammatic method we present here is  extensions
of quandles.
Cohomology theories of groups and other algebraic systems
have interpretations in terms of group extensions or
obstructions to deformations of algebraic systems
(see for example \cite{Brown,GS} and for a diagrammatic approach
\cite{MS}).
We give analogous interpretations for quandle cohomology.
The proofs are based on knot diagrams.

The paper is organized as follows.
In Section~1, we give a summary of preliminary material
  on quandle homology and cocycle knot invariants.
Applications of the Kronecker product are given in Section~2,
and the extensions of quandles by cocycles are investigated in Section~3.

{\bf Acknowledgement.\/} We are grateful to Edwin Clark
for  his valuable comments.

\section{Quandle homology and colored knot diagrams}

In this section we review necessary material from the papers mentioned
in the introduction.

A {\it quandle}, $X$, is a set with a binary operation $(a, b) \mapsto a * b$
such that

(I) For any $a \in X$,
$a* a =a$.

(II) For any $a,b \in X$, there is a unique $c \in X$ such that
$a= c*b$.

(III)
For any $a,b,c \in X$, we have
$ (a*b)*c=(a*c)*(b*c). $

A {\it rack} is a set with a binary operation that satisfies
(II) and (III).
Racks and quandles have been studied in, for example,
\cite{Brieskorn,FR,Joyce,K&P,Matveev}.
The axioms for a quandle correspond respectively to the
Reidemeister moves of type I, II, and III
(see
\cite{FR}, \cite{K&P}, for example).

A function $f: X \rightarrow  Y$ between quandles
or racks  is a {\it homomorphism}
if $f(a \ast b) = f(a) * f(b)$
for any $a, b \in X$.

The following are typical examples of quandles.

\begin{itemize}
\item
A group $X=G$ with
$n$-fold conjugation
as the quandle operation: $a*b=b^{-n} a b^n$.
\item
Any set $X$ with the operation $x*y=x$ for any $x,y \in X$ is
a quandle called the {\it trivial} quandle.
The trivial quandle of $n$ elements is denoted by $T_n$.
\item
Let $n$ be a positive integer.
For elements  $i, j \in \{ 0, 1, \ldots , n-1 \}$, define
$i\ast j \equiv 2j-i \pmod{n}$.
Then $\ast$ defines a quandle
structure  called the {\it dihedral quandle},
  $R_n$.
This set can be identified with  the
set of reflections of a regular $n$-gon
  with conjugation
as the quandle operation.
\item
Any $\Lambda (={\bf Z}[T, T^{-1}])$-module $M$
is a quandle with
$a*b=Ta+(1-T)b$, $a,b \in M$, called an {\it  Alexander  quandle}.
Furthermore for a positive integer
$n$, a {\it mod-$n$ Alexander  quandle}
${\bf Z}_n[T, T^{-1}]/(h(T))$
is a quandle
for
a Laurent polynomial $h(T)$.
The mod-$n$ Alexander quandle is finite
if the coefficients of the
highest and lowest degree terms
of $h$
  are $\pm 1$.
\end{itemize}

See \cite{Brieskorn,FR,Joyce,Matveev}
  for further examples of quandles.

\begin{figure}
\begin{center}
\mbox{
\epsfxsize=3in
\epsfbox{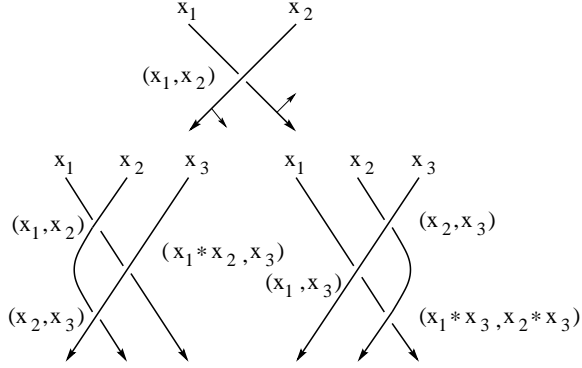}
}
\end{center}
\caption{ Type III move and the quandle identity  }
\label{nu2cocy}
\end{figure}

\bigskip

  Let $C_n^{\rm R}(X)$ be the free
abelian group generated by
$n$-tuples $(x_1, \dots, x_n)$ of elements of a quandle $X$. Define a
homomorphism
$\partial_{n}: C_{n}^{\rm R}(X) \to C_{n-1}^{\rm R}(X)$ by \begin{eqnarray}
\lefteqn{
\partial_{n}(x_1, x_2, \dots, x_n) } \nonumber \\ && =
\sum_{i=2}^{n} (-1)^{i}\left[ (x_1, x_2, \dots, x_{i-1}, x_{i+1},\dots, 
x_n) \right.
\nonumber \\
&&
- \left. (x_1 \ast x_i, x_2 \ast x_i, \dots, x_{i-1}\ast x_i, x_{i+1}, 
\dots, x_n) \right]
\end{eqnarray}
for $n \geq 2$
and $\partial_n=0$ for
$n \leq 1$.
  Then
$C_\ast^{\rm R}(X)
= \{C_n^{\rm R}(X), \partial_n \}$ is a chain complex.

Let $C_n^{\rm D}(X)$ be the subset of $C_n^{\rm R}(X)$ generated
by $n$-tuples $(x_1, \dots, x_n)$
with $x_{i}=x_{i+1}$ for some $i \in \{1, \dots,n-1\}$ if $n \geq 2$;
otherwise let $C_n^{\rm D}(X)=0$. If $X$ is a quandle, then
$\partial_n(C_n^{\rm D}(X)) \subset C_{n-1}^{\rm D}(X)$ and
$C_\ast^{\rm D}(X) = \{ C_n^{\rm D}(X), \partial_n \}$ is a sub-complex of
$C_\ast^{\rm
R}(X)$. Put $C_n^{\rm Q}(X) = C_n^{\rm R}(X)/ C_n^{\rm D}(X)$ and
$C_\ast^{\rm Q}(X) = \{ C_n^{\rm Q}(X), \partial'_n \}$,
where $\partial'_n$ is the induced homomorphism.
Henceforth, all boundary maps will be denoted by $\partial_n$.

For an abelian group $G$, define the chain and cochain complexes
\begin{eqnarray}
C_\ast^{\rm W}(X;G) = C_\ast^{\rm W}(X) \otimes G, \quad && \partial =
\partial \otimes {\rm id}; \\ C^\ast_{\rm W}(X;G) = {\rm Hom}(C_\ast^{\rm
W}(X), G), \quad
&& \delta= {\rm Hom}(\partial, {\rm id})
\end{eqnarray}
in the usual way, where ${\rm W}$
  $={\rm D}$, ${\rm R}$, ${\rm Q}$.

The $n$\/th {\it quandle homology group\/}  and the $n$\/th
{\it quandle cohomology group\/ } \cite{CJKLS} of a quandle $X$ with 
coefficient group $G$ are
\begin{eqnarray}
H_n^{\rm Q}(X; G)
  = H_{n}(C_\ast^{\rm Q}(X;G)), \quad
H^n_{\rm Q}(X; G)
  = H^{n}(C^\ast_{\rm Q}(X;G)). \end{eqnarray}

\begin{sloppypar}
The cycle and boundary groups
(resp. cocycle and coboundary groups)
are denoted by $Z_n^{\rm Q}(X;G)$ and $B_n^{\rm Q}(X;G)$
(resp.  $Z^n_{\rm Q}(X;G)$ and $B^n_{\rm Q}(X;G)$),
  so that
$$H_n^{\rm Q}(X;G) = Z_n^{\rm Q}(X;G)/ B_n^{\rm Q}(X;G),
\; H^n_{\rm Q}(X;G) = Z^n_{\rm Q}(X;G)/ B^n_{\rm Q}(X;G).$$
We will omit the coefficient group $G$  as usual if $G = {\bf Z}$.
\end{sloppypar}

\bigskip

Let a classical knot diagram
be given.
The co-orientation is a family of normal vectors to the knot diagram
such that the pair (orientation, co-orientation) matches
the given (right-handed, or counterclockwise) orientation of the plane.
At a crossing,
if the pair of the co-orientation
  of the
over-arc and  that of the under-arc
matches the (right-hand) orientation of the plane, then the
crossing is called {\it positive}; otherwise it is {\it negative}.
The crossings  depicted
in Fig.~\ref{nu2cocy} are positive by convention.

A  {\it coloring}
of an oriented  classical knot diagram is a
function ${\mathcal C} : R \rightarrow X$, where $X$ is a fixed
quandle
and $R$ is the set of over-arcs in the diagram,
satisfying the  condition
depicted
in the top
of Fig.~\ref{nu2cocy}.
In the figure, a
crossing with
over-arc, $r$, has color ${\mathcal C}(r)= y \in X$.
The under-arcs are called $r_1$ and $r_2$ from top to bottom;
the normal (co-orientation) of the over-arc $r$ points from $r_1$ to $r_2$.
Then it is required that
${\mathcal C}(r_1)= x$ and ${\mathcal C}(r_2)=x*y$.

Note that locally the colors do not depend on the
orientation of the under-arc.
The quandle element ${\mathcal C}(r)$ assigned to an arc $r$ by a coloring
  ${\mathcal C}$ is called a {\it color} of the arc.
This definition of colorings on knot diagrams has been known, see
\cite{FR,FoxTrip} for example.
Henceforth, all the quandles that are used to color diagrams will be finite.

In Fig.~\ref{nu2cocy} bottom, the relation between
the
Redemeister type III move
and
quandle axiom (self-distributivity) is indicated.
In particular, the colors of the bottom right segments before and after
the move correspond to
self-distributivity.

\begin{figure}
\begin{center}
\mbox{
\epsfxsize=3in
\epsfbox{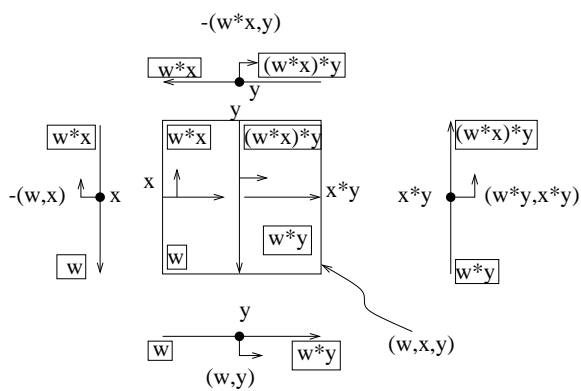}
}
\end{center}
\caption{Representing a $3$-chain by shadow colored diagrams }
\label{3chain}
\end{figure}

A {\it shadow coloring} (or {\it face coloring}) of
a classical knot diagram is a
function ${\mathcal C} : \tilde{R} \rightarrow X$, where $X$ is a fixed
quandle
and $ \tilde{R}$ is the set of arcs in the diagram and regions separated by
the underlying immersed curve of the knot diagram,
satisfying the  condition
depicted
in the middle square
of Fig.~\ref{3chain}.
In the figure,
arcs are colored under the same rule as above, and the regions are also
colored by the following similar rule.
Let $R_1$ and $R_2$ be the regions separated by an arc $r$ colored by $x$.
Suppose that the normal to $r$ points from $R_1$ to $R_2$.
If $R_1$ is colored by $w$, then $R_2$ is required to be colored by $w*x$.
Note that near a crossing  there are more than one way
  to go from one region to
another, but the
self-distributivity guarantees
unique colors near a crossing.
In the figures,  colors on the the regions are depicted as
letters enclosed within squares.

\begin{figure}
\begin{center}
\mbox{
\epsfxsize=3in
\epsfbox{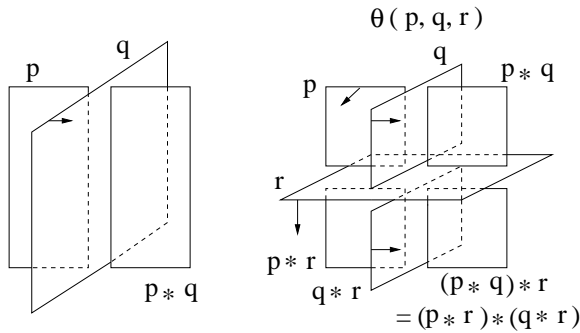}
}
\end{center}
\caption{Colors at double curves and triple points }
\label{triplepoint}
\end{figure}

Colorings and shadow colorings are defined for knotted surfaces
in $4$-space similarly using their diagrams in $3$-space.
The coloring rule is depicted in Fig.~\ref{triplepoint}.

\bigskip

\begin{figure}
\begin{center}
\mbox{
\epsfxsize=2in
\epsfbox{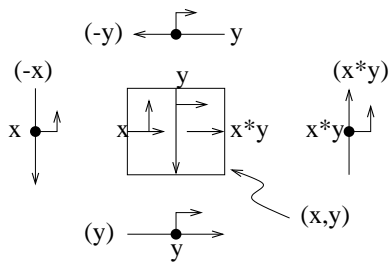}
}
\end{center}
\caption{Representing a $2$-chain by  colored diagrams }
\label{2chain}
\end{figure}

Each positive
crossing in a colored knot diagram represents a
pair $(x,y) \in C_2^{\rm R}(X;A)$ as depicted in Fig.~\ref{nu2cocy} top.
The first factor $x$ is the color on an under-arc away from which
the normal of the over-arc points. The color $y$ is on the over-arc.
If the crossing is negative, then such colors represent  $-(x,y)$.
The relation between colored crossings and boundary homomorphisms
is depicted in Fig.~\ref{2chain}.
The $1$-chains
  represented by colored and oriented points on
line segments are depicted as boundaries in Fig.~\ref{2chain}.
The color  on the center vertex on an arc determine the $1$-chain
that the vertex represents.
These colors are the colors at the end points of the
arc in the crossing. The sign of the $1$-chain is determined by pushing
the normal to the arc into the boundary and comparing to the oriented
subarc of the boundary given the counterclockwise orientation.
In Fig.~\ref{2chain}, the boundary terms of the $2$-chain $(x,y)$ are
$1$-chains $(x*y)$,  $(-y)$,  $(-x)$, and $(y)$, and thier formal sum
matches the negative of  $\partial (x,y)$, where
$\partial$ is the boundary homomorphism of quandle homology.
In particular, any  colored knot diagram represents a $2$-cycle,
as the boundary terms cancel.

Similarly, shadow colored crossings represent triples
$(w,x,y)$ as depicted in Fig.~\ref{3chain}.
The boundaries are also indicated in the figure.
In particular, shadow colored knot diagrams represent
$3$-cycles.
The signs are determined as above and the chain (for example $-(w,x)$ on
the left) is determined as follows. The color $w$ is the color in the
region
away from which the normal to the
arc colored $x$ points,
and $x$ is the color on that arc.

\begin{figure}
\begin{center}
\mbox{
\epsfxsize=1.2in
\epsfbox{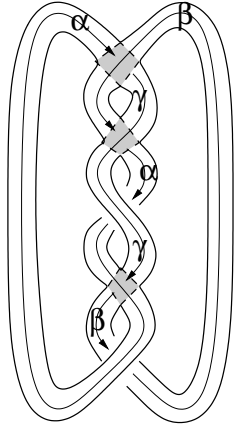}
}
\end{center}
\caption{A $1$-knot diagram }
\label{ssvsv1}
\end{figure}

Colored or shadow colored classical knot diagrams in orientable surfaces
are defined similarly, and represent $2$- or $3$-cycles, respectively,
as boundary terms also cancel.
For example, the colored knot diagram on a surface depicted
  in Fig.~\ref{ssvsv1}
represents a $2$-cycle
$c=(\alpha, \beta) + (\beta, \gamma)- (\beta, \alpha)$
in $Z^2_{\rm Q}(R_3; {\bf Z}_3)$,
  where $\{\alpha, \beta, \gamma \}=\{ 0,1,2\} = R_3$.
In fact, $[c]=0 \in H^2_{\rm Q}(R_3; {\bf Z}_3)$, as it is known
\cite{CJKLS} that
$H^2_{\rm Q}(R_3; {\bf Z}_3)=0$.
The shaded regions in Fig.~\ref{ssvsv1} are crossings
of the diagram, and the unshaded
crossing in the middle is a cross-over of the surface, and is not
a crossing of the diagram.
Since all boundaries of each colored crossing are attached to other
colored crossings, we see that the boundary terms cancel,
  and the diagram represents
a $2$-cycle.

\begin{figure}
\begin{center}
\mbox{
\epsfxsize=.7in
\epsfbox{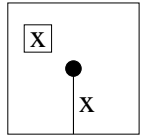}
}
\end{center}
\caption{The endpoint diagram with colors}
\label{endpoint}
\end{figure}

\begin{figure}
\begin{center}
\mbox{
\epsfxsize=4in
\epsfbox{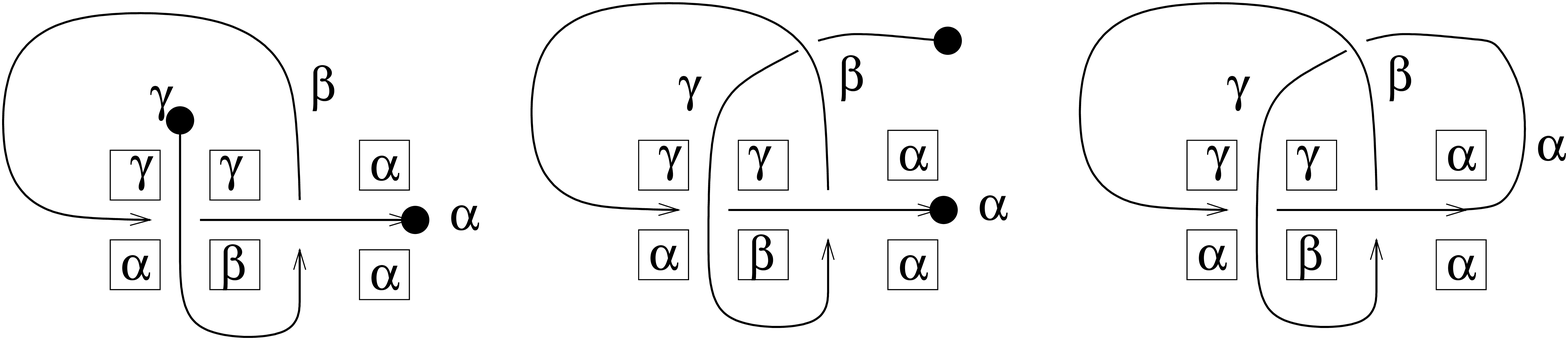}
}
\end{center}
\caption{Generators of $H^3_{\rm Q}(R_3;{\bf Z}_3)$  }
\label{gentre}
\end{figure}

Furthermore, for shadow colored diagrams, we allow the colored end point
diagram that is
depicted in Fig.~\ref{endpoint}.
Since the boundary term of the end point diagram represents $\pm (x,x)$,
which represents zero as a quandle chain, a shadow colored knot diagram
represents a quandle $3$-cycle even with such colored end points allowed.
Examples are  depicted in Fig.~\ref{gentre} left and middle.
All three diagrams in Fig.~\ref{gentre}
  represent a generator of  $H^3_{\rm Q}(R_3;{\bf Z}_3)$.

These colored knot diagrams on surfaces
(possibly with end points for shadow colors)
  are called  colored {\it abstract} knot (or arc, respectively)
diagrams.

\begin{figure}
\begin{center}
\mbox{
\epsfxsize=3in
\epsfbox{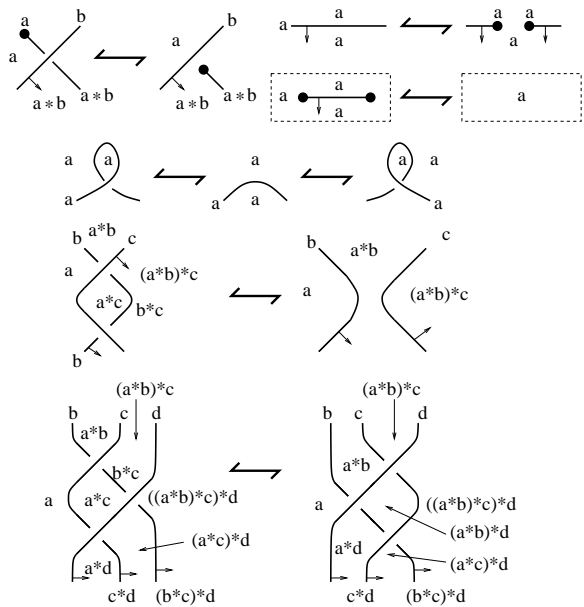}
}
\end{center}
\caption{Moves for shadow colored diagrams, Part I (Reidemeister moves)}
\label{shmoves2}
\end{figure}

\begin{figure}
\begin{center}
\mbox{
\epsfxsize=4.2in
\epsfbox{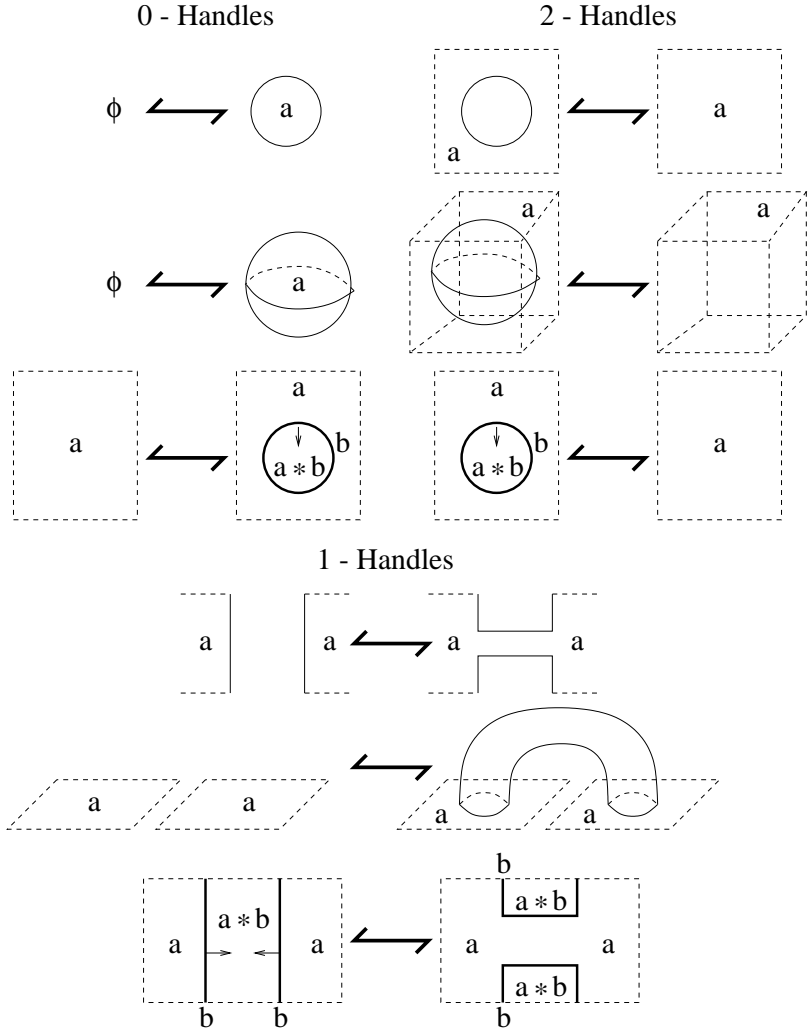}
}
\end{center}
\caption{Moves for shadow colored diagrams, Part II (Morse critical points) }
\label{shmoves1}
\end{figure}

Colored and shadow colored abstract knot diagrams
represent cycles, and the moves for colored diagrams are known \cite{SSS2}:
Colored diagrams related by a finite sequence of the moves
represent the same homology class.
  Such moves are depicted in Figs.~\ref{shmoves2} and \ref{shmoves1},
  and called the quandle homology moves.

\bigskip

The (fundamental)
quandle of an $n$-knot
diagram \cite{Joyce,Matveev} is generated by the $n$-regions of the diagram;
  the relations in the quandle can be
read  from the crossings ($n=1$) or double point curves ($n=2$).
See \cite{FR,K&P} for Wirtinger presentations
of knot quandles defined from
knot diagrams, which
are
similar to Wirtinger presentations of knot groups.
In this case, arcs of knot diagrams represent generators, and crossings
give relations of Wirtinger form.
Let
$Q(K)$
represent such a quandle
for a knot diagram
$K.$

Let $K$ be a knot diagram on a compact oriented surface $F$.
Then the {\it fundamental
shadow quandle} $SQ(K)$ is defined as follows \cite{SSS2}.
The generators correspond to
over-arcs and connected components
of $F \setminus ({\mbox{\rm universe of}} \ K)$,
where the universe is the underlying immersed curves of $K$
(without crossing information).
The relations are defined for
each crossing as ordinary
fundamental quandles, and at each arc dividing
regions. Specifically, if $a$ and
$b$ are generators corresponding to
adjacent regions such that the
normal points from the region colored $a$
to that colored $b$, and if the arc dividing these
regions is colored by $c$, then we have the relation $b=a*c$.
This defines a presentation of a quandle, which is called the
fundamental shadow quandle of $K$.
Two diagrams on $F$ that differ by Reidemeister moves on $F$
have
isomorphic fundamental shadow quandles.
The shadow colors are regarded as  quandle homomorphisms
from the fundamental shadow quandle to a
quandle $X$.

\section{Cocycle invariants as Kronecker products}

In this section we give an interpretation of
the state-sum invariants (called
quandle cocycle invariants) defined in \cite{CJKLS}
in terms of pairings on quandle homology theory.
Using this interpretation, a new method of computing these invariants
is given.
In this section all (co)chain, (co)cycle, (co)boundary,
  and (co)homology groups
are quandle groups and we
sometimes
drop the letter ${\rm Q}$ from the
subscript or superscript.



Let $({\mathcal C}, \partial)$ be a
  chain complex with the boundary homomorphism $\partial$.
Let $\langle \quad , \quad \rangle : Z_n({\mathcal C} ) \otimes 
Z^n({\mathcal C};G)$
be the Kronecker product, where $G$ is a coefficient abelian group,
and it is omitted
as usual if $G={\bf Z}$.
Thus  $\langle \eta , \phi \rangle = \phi(\eta)$
for any $\eta \in  Z_n({\mathcal C} )$ and $\phi \in Z^n({\mathcal C} ;G) $.
  This pairing induces a well-defined bilinear pairing (Kronecker product)
$\langle \quad , \quad \rangle : H_n({\mathcal C} ) \otimes H^n({\mathcal 
C} ;G)
\rightarrow G$.

Let $({\mathcal C}_0, \partial)$ be another chain complex,
and denote by $\mbox{Hom}({\mathcal C}_0, {\mathcal C})$
be the set of chain homomorphisms.
For $\eta \in Z_n( {\mathcal C}_0)$ and $\phi \in Z^n( {\mathcal C} ; G)$,
define
$$\Phi(\eta, \phi)
  = \sum_{ \{ f \in {\mathcal F} \subset \mbox{Hom}({\mathcal C}_0, 
{\mathcal C}) \} }
\langle  f_{\#}\eta ,  \phi \rangle ,$$
where ${\mathcal F}$ is a fixed finite
subset of $\mbox{Hom}({\mathcal C}_0, 
{\mathcal C})$.
This defines a bilinear pairing
$\Phi: Z_n({\mathcal C}_0) \otimes Z^n({\mathcal C};G) \rightarrow {\bf Z}[G]$.
Since each Kronecker product depends only
  on the homology and cohomology classes, we have the following.

\begin{lemma}
The above defined $\Phi$ does not depend on the choice of (co)cycles
and is determined only by their (co)homology classes. Thus it induces
a well-defined bilinear pairing
$$ \Phi([\eta], [\phi]) = \sum_{ \{ f \in  {\mathcal F} \subset 
\mbox{Hom}({\mathcal C}_0, {\mathcal C})  \} }
\langle f_{*} [\eta] ,  [\phi] \rangle .$$
\end{lemma}

\begin{definition}
Let $K$ be an abstract
  $n$-knot diagram, for $n=1,2$.
In the above description, let ${\mathcal C}$ be the chain complex
$\{ C_*^{\rm Q}(X) \}$ for a finite quandle $X$,
and  ${\mathcal C}_0$ be $\{ C_*^{\rm Q}(\Pi) \}$ where
$\Pi=Q(K)$ is the fundamental quandle of a knot $K$.
A knot diagram $K$ represents a class  $[K] \in H_{n+1}(\Pi)$.
Pick and fix $\phi \in Z^{n+1}(X; G)$.
Let ${\mathcal F} \subset \mbox{Hom}({\mathcal C}_0, {\mathcal C})$ be
the set of all  chain maps induced from all quandle homomorphisms
$\Pi \rightarrow X$,
that is,
$${\mathcal F}= \{ \ f_{\#} : {\mathcal C}_0  \rightarrow {\mathcal C}\  |\
f: \Pi \rightarrow X : \mbox{quandle homomorphism} \  \} . $$
Define
$ \Phi_{\phi}(K)= \Phi ([K], [\phi]) \in {\bf Z}[G] $.
This is called the {\it quandle cocycle invariant} of $K$ with
color quandle $X$.

A similar invariant,
called {\it shadow quandle cocycle invariant},
  is defined using fundamental shadow quandles $SQ(K)$,
by $S\Phi_{\theta}(K)=\Phi  ([K], [\theta]) \in {\bf Z}[G]$,
where $[K] \in H_{n+2}(SQ(K))$ and $\theta \in  H^{n+2}(X; G)$.
  \end{definition}

\begin{remark}
The (shadow) quandle cocycle invariants coincide
with the state-sum invariant defined in \cite{CJKLS}.
The set of colorings in \cite{CJKLS} corresponds to
$ {\mathcal F} $, and the Boltzmann weight defined
from a fixed cocycle $\phi$ corresponds to
$\langle [K], f^{\#} \phi \rangle$.
The definition in terms of pairing was suggested to us by Fenn and Rourke
in a correspondence.
The above gives a generalization using quandle
homology and abstract knots.
The following generalizes the invariants to abstract knots.
  \end{remark}

\begin{proposition}
The quandle cocycle invariant $ \Phi_{\phi}(K)$ is an invariant
for abstract knots of dimensions $3$ and $4$
(i.e., it does not depend on the choice of the diagram,
and is well-defined up to equivalence of abstract knot diagrams).

The shadow  quandle cocycle invariant $S \Phi_{\phi}(K)$ is an invariant
of knot diagrams in Euclidean spaces ${\bf R}^n$ up to
Reidemeister moves (and their analogues in dimension $4$).
\end{proposition}
\begin{proof}
The equivalence relation in question does not alter the class
$[K] \in H_n(Q(K))$ or $ H_{n+1}(SQ(K))$ in the situations stated.
\end{proof}

The above interpretation can be used for computation of
the invariant as follows.

Recall that $R_3$ is ${\bf Z}_3$ as a set with the quandle operation
$a*b \equiv 2b-a$ (mod 3).
Let $\Phi(K)$ denote the cocycle invariant of classical knots
defined by shadow colorings by $R_3$ and the cocycle
  $\xi \in Z^3_{\rm Q} (R_3; {\bf Z}_3)$ defined by
$$\xi = \chi_{012} \chi_{021}
\chi_{101} \chi_{201} \chi_{202}\chi_{102}$$
where
$$\chi_{abc} (x,y,z) = \left\{ \begin{array}{lr} t & {\mbox{\rm if }} \
(x,y,z)=(a,b,c), \\
1 & {\mbox{\rm if }} \
(x,y,z)\not=(a,b,c). \end{array}\right.$$
In this section we give examples of this invariant using the pairing
interpretations given in the previous section.

\begin{lemma} \label{allcolorlem}
For any $\{ \alpha, \beta, \gamma \} = \{ 0,1,2\}$,
the shadow colored diagram
on the left of
Fig.~\ref{gentre}
represents a generator of $H_3^{\rm Q}(R_3;{\bf Z}_3) \cong {\bf Z}_3 $.

For any $\{ \alpha, \beta, \gamma \} = \{ 0,1,2\}$,
the shadow colored diagram which is the mirror image of that in
the left of
Fig.~\ref{gentre}
represents the inverse of the  generator of
  $H_3^{\rm Q}(R_3;{\bf Z}_3) \cong {\bf Z}_3 $ given above.
\end{lemma}
\begin{proof}
The shadow colored diagram represents
$(\alpha, \beta, \gamma) + (\alpha, \gamma, \alpha)$.
It is known \cite{CJKLS} that
the cocycle
$\xi$
which is defined above
represents a generator of  $H^3_{\rm Q}(R_3;{\bf Z}_3) \cong {\bf Z}_3 $.
All possibilities of $\{ \alpha, \beta, \gamma \} = \{ 0,1,2\}$
are
evaluated by
$\xi$
to give the
generator
$1 \in {\bf Z}_3$,
and the result follows. The mirror image is checked similarly.
\end{proof}

Note that
the
other diagrams in Fig.~\ref{gentre}
also
represent
the same generator, and
have
the same property as stated in the above lemma.

\begin{figure}
\begin{center}
\mbox{
\epsfxsize=2in
\epsfbox{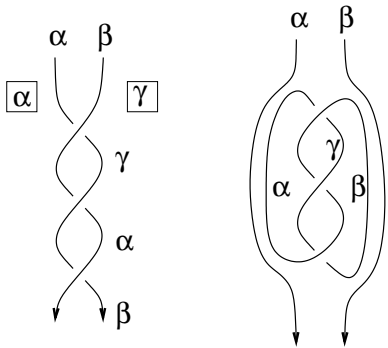}
}
\end{center}
\caption{A decomposition of $T(2, n)$ }
\label{torus2}
\end{figure}

\begin{proposition}
A torus knot or  link $T(2, n)$ is $3$-colorable if and only if
$n$ is a multiple of $3$, $n=3k$ for some integer $k$.
In this case,
$\Phi( T(2, 3k) ) = 9 + 18 t^k $. Here $k$ is regarded as an element
of ${\bf Z}_3$.
\end{proposition}
\begin{proof}
We represent $T(2,n)$ as a closed $2$-braid as depicted in Fig.~\ref{torus2}.
If the top two segments receive the same color, the cocycle
invariant is trivial. There are $9$ such shadow colors.
If the two top segments receive distinct colors, then perform the
quandle homology moves as indicated in the figure, and produce a copy of
a generator in Fig.~\ref{gentre} for
each set of  three crossings. By Lemma~\ref{allcolorlem},
for any choice of $\{ \alpha, \beta, \gamma \}= \{ 0,1,2\}$,
the figure represents the generator of $H_3(R_3 ; {\bf Z}_3)$.
Hence for any non-trivial coloring ${\mathcal C}$,
$\langle  [T(2, 3k)], {\mathcal C}^* \xi \rangle = t^k$.
The result follows.
\end{proof}

\begin{figure}
\begin{center}
\mbox{
\epsfxsize=2.5in
\epsfbox{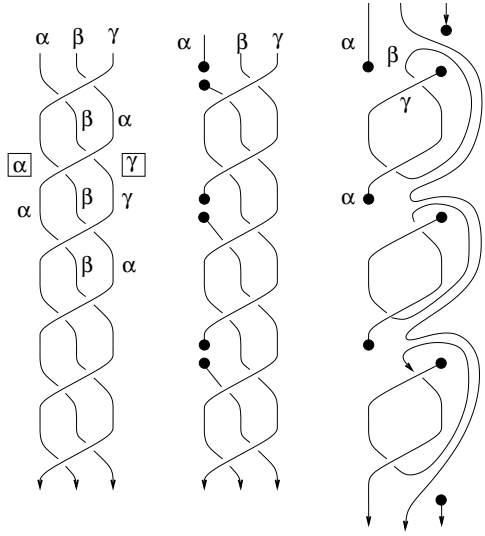}
}
\end{center}
\caption{A decomposition of $T(3,n)$  }
\label{torus31}
\end{figure}

\begin{figure}
\begin{center}
\mbox{
\epsfxsize=3.5in
\epsfbox{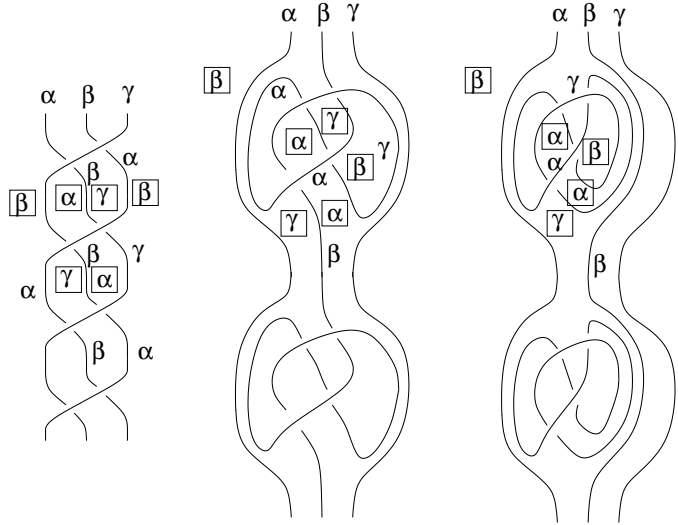}
}
\end{center}
\caption{Another decomposition of $T(3,n)$  }
\label{torus32}
\end{figure}

\begin{figure}
\begin{center}
\mbox{
\epsfxsize=3.5in
\epsfbox{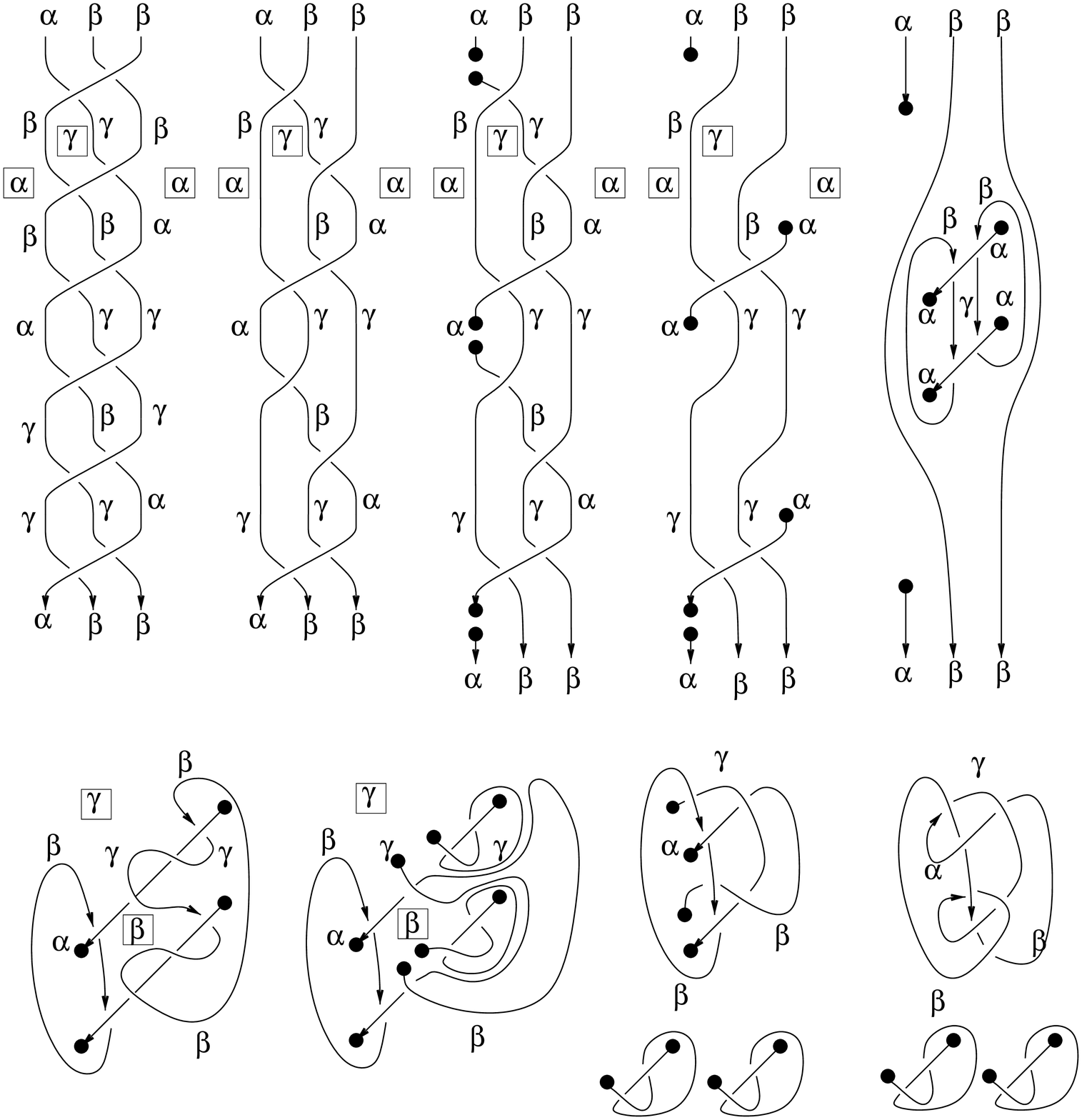}
}
\end{center}
\caption{A coloring with repetitive colors for  $T(3,n)$, type I  }
\label{torus33}
\end{figure}

\begin{figure}
\begin{center}
\mbox{
\epsfxsize=3.5in
\epsfbox{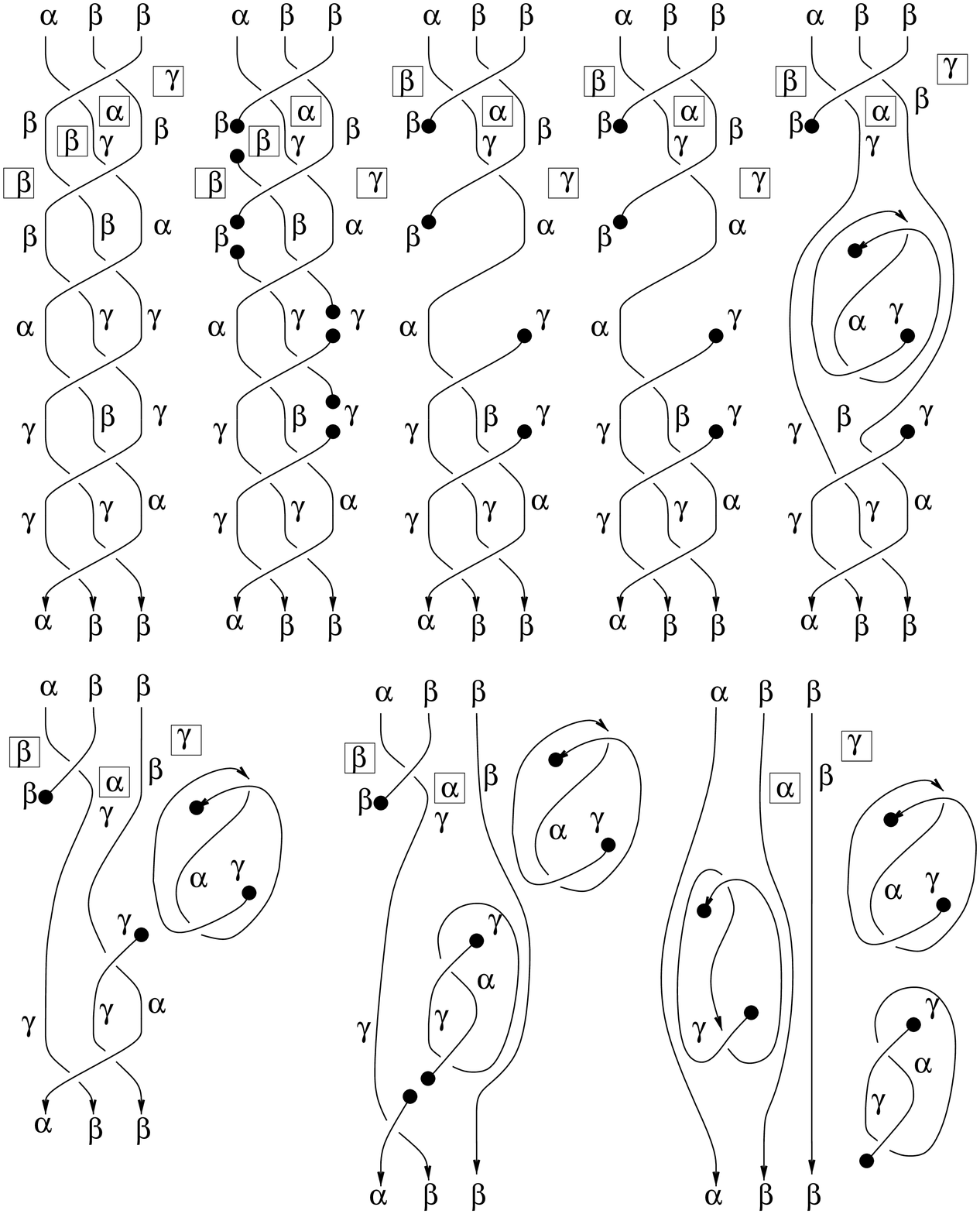}
}
\end{center}
\caption{A coloring with repetitive colors for  $T(3,n)$, type II }
\label{torus34}
\end{figure}

\begin{figure}
\begin{center}
\mbox{
\epsfxsize=3.5in
\epsfbox{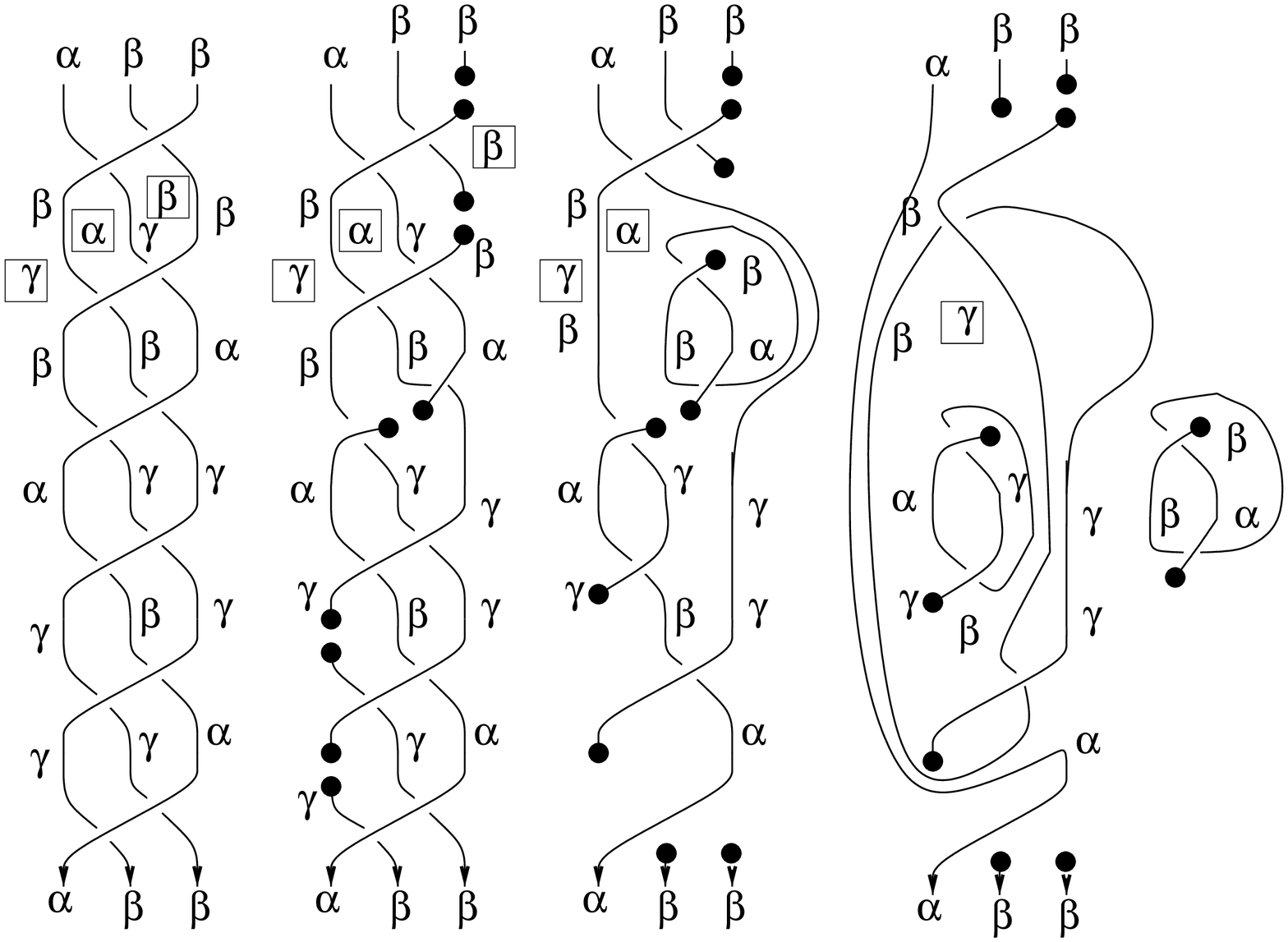}
}
\end{center}
\caption{A coloring with repetitive colors for  $T(3,n)$, type III }
\label{torus35}
\end{figure}

\begin{proposition}
A torus knot or link $T(3, n)$ is $3$-colorable if and only if
$n$ is a multiple of $2$, $n=2k$ for some integer $k$.
If $k$ is not a multiple of $3$, then
$\Phi( T(3, 2k) ) = 9 +  18 t^k $,
where  $k$ is regarded as an element
of ${\bf Z}_3$. Otherwise, $\Phi( T(3, 2k) ) =45$.
\end{proposition}
\begin{proof}
There are three cases to color the top strings.
All three  colors are distinct, two colors are the same, and all
colors
are the same.
If all colors are the same, such colorings contribute $1$ to the invariant,
and there are $9$ such shadow colorings.
If all colors are distinct, then such a coloring is depicted
  in Fig.~\ref{torus31} left. Such a color exists
if and only if
$n$
is even.
Note that the color of the middle string,
$\beta$, stays in the middle strings. Hence we need to consider two types
of shadow colors, (1) $\beta$ is not a color of the unbounded region,
(2) it is.
The case (1), (2)  are  shown in Fig.~\ref{torus31}, \ref{torus32}
respectively. For each case, the figure shows that such a colored
diagram is homologous to a copy of a generator for each set of $6$ crossings.
Hence this contributes $t^k$ to the invariant.

By symmetry of the diagram, if two colors are the same, it can be assumed that
the top strings receive, say, $(\alpha, \beta, \beta)$
in this order from left
to right,
as depicted in
the left hand side of
Fig.~\ref{torus33}.
In this case $k$ must be a multiple of
$3$ for such a color to exist, as can be seen from the figure.
Hence the first case (when $k$ is not a multiple of $3$) follows from the
above argument, and the cases
for  $T(3, 3\ell)$ remain.  A block of one contribution of $\ell$
(i.e., $T(3,3)$)
is depicted in the figure.
There are $6$ such colorings, and there are $18$ such shadow colorings.
There are three cases for shadow colors, as depicted in Figs.~\ref{torus33},
\ref{torus34}, \ref{torus35} left, respectively.
For each of these cases, the figures show that they are homologous to
three copies of a generator, contributing $1$ to the invariant
(as each contribution is counted modulo $3$).
The second case follows.
\end{proof}

\begin{figure}
\begin{center}
\mbox{
\epsfxsize=2.5in
\epsfbox{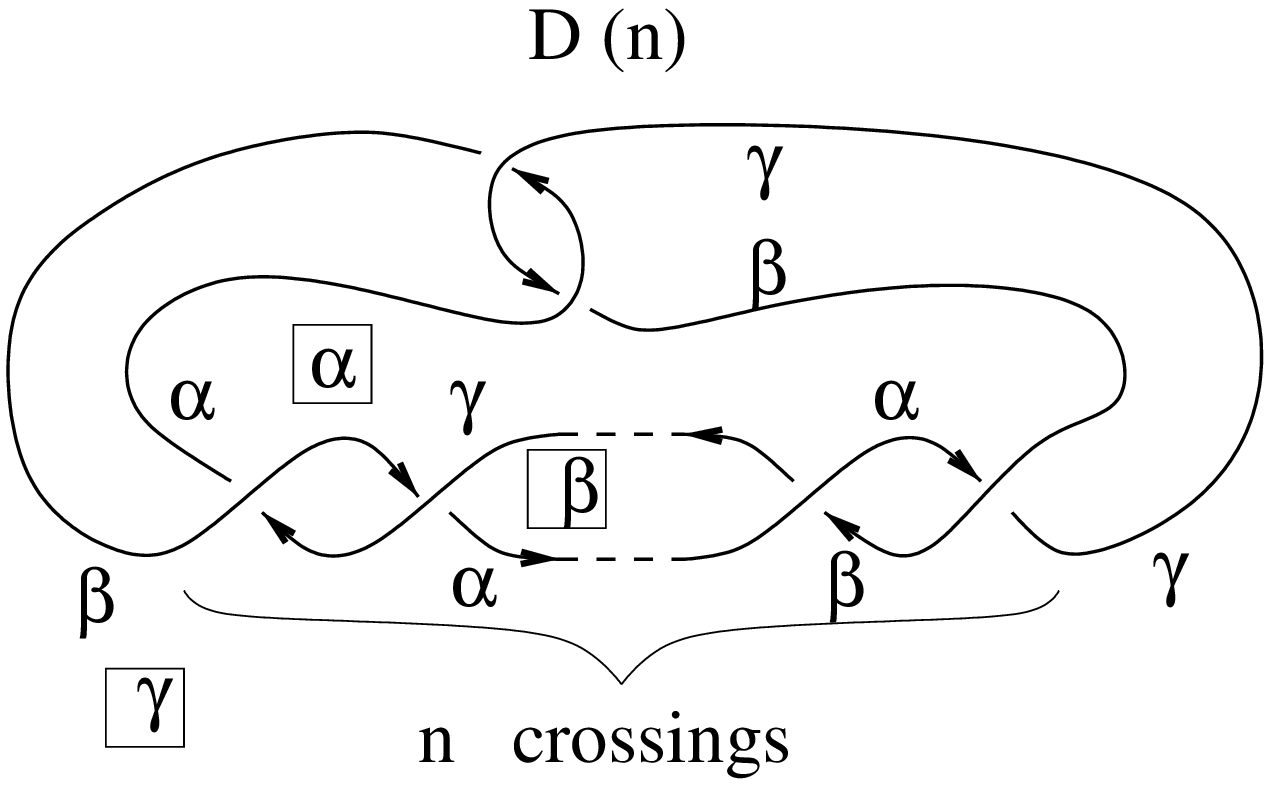}
}
\end{center}
\caption{Doubled knots $D(n)$}
\label{dble}
\end{figure}

Next we compute the invariant for doubled knots $D(n)$
depicted in Fig.~\ref{dble}.
These knots are twisted  Whitehead doubles of the unknot.
The integer $n$ represents the number of crossings as indicated.
In the figure the crossings are positive ones, and if $n$ is
negative, we take negative crossings.

\begin{figure}
\begin{center}
\mbox{
\epsfxsize=3in
\epsfbox{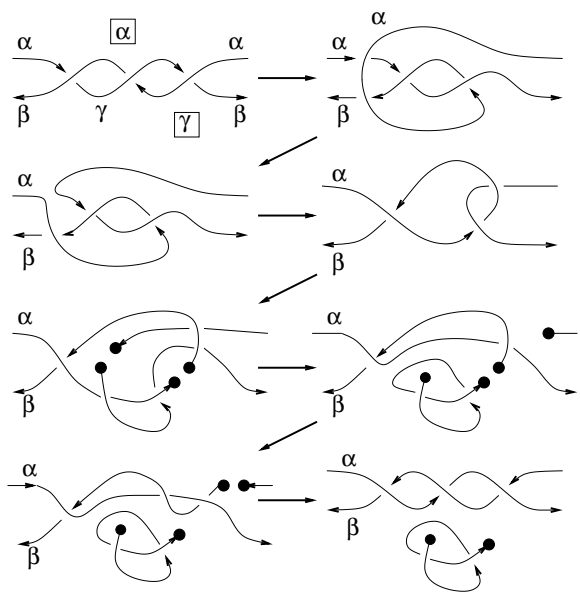}
}
\end{center}
\caption{Decomposing antiparallel strings}
\label{antipara}
\end{figure}

\begin{proposition}
A doubled knot $D(n)$ is $3$-colorable if and only if
$n=3k+1$. If $k=2m$, then   $\Phi(D(3k+1))= 9 + 18 t^{m+1}$,
and if $k=2m+1$, then  $\Phi(D(3k+1))= 9 + 18 t^{m}$.
\end{proposition}
\begin{proof}
Figure~\ref{dble} shows the first half of the statement.
Note that there are $3$ trivial colorings on strings and $9$ corresponding
shadow colorings, and $18$ shadow colorings corresponding to non-trivial
colorings on strings.
In Fig.~\ref{antipara}, it is shown that for a particular shadow coloring,
a set of three crossings in the diagrams of $D(n)$ contributes
a single copy of a generator of $H^3_{\rm Q}(R_3;{\bf Z}_3)$,
and changes the crossings from positive to negative. These three negative
crossings cancel with the next set of three positives.
Hence this particular shadow coloring contributes $t^{m}$ from
$6m$ crossings. If $n=6m+1$, then
the rest is a single crossing, and $D(1)$ is a right-handed trefoil knot,
giving another generator.
For other choices of colorings on the regions, simply change the starting
point of replacement and use the same argument.
Apply type II Reidemeister moves if necessary, to reduce
the given colored diagram to $n=1$ case.
This gives the case $n=6m+1$, $\Phi(D(n))=9 + 18 t^{m+1}$.
If $n=6m+4$, then apply the operation in Fig.~\ref{antipara} $(m+1)$ times,
to obtain $D(-2)$, which is a left-handed trefoil, representing the negative
of the generator. Hence the contribution is $(m+1)-1=m$, giving
  $\Phi(D(6m+4))= 9 + 18 t^{m}$.
\end{proof}

\begin{figure}
\begin{center}
\mbox{
\epsfxsize=3in
\epsfbox{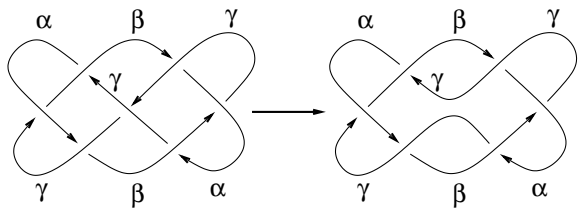}
}
\end{center}
\caption{Deforming $7_4$}
\label{74}
\end{figure}

Let $T'(2n)$ denote the type $(2, 2n)$ torus link of $2$-components
with opposite orientations given to parallel strings.
In other words, $T'(2n)$ is obtained from Fig.~\ref{antipara}
with $2n$ crossings by taking the ``braid closure'' between
left and right ends.

\begin{corollary}
The link $T'(2n)$ is $3$-colorable if and only if $n=3k$.
In this case, $\Phi(T'(6k))=9 + 18 t^k$.
\end{corollary}
\begin{proof}
The method depicted in Fig.~\ref{antipara} applies
in the same manner as in the proof of the above proposition.
\end{proof}

\begin{remark}
The methods developed  so far can be applied effectively
to other examples to evaluate the invariant,
directly or indirectly.
  As examples, we examine knots in the table.
There are four knots in the table less than $8$ crossing that
are $3$-colorable: $3_1$, $6_1$,  $7_4$, and $7_7$.

\begin{figure}
\begin{center}
\mbox{
\epsfxsize=2.5in
\epsfbox{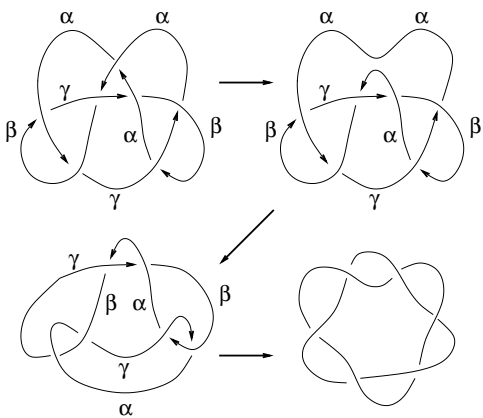}
}
\end{center}
\caption{Deforming  $7_7$}
\label{77}
\end{figure}

\begin{itemize}
\item
$3_1$ is a generator, so $\Phi(3_1)=9 + 18t$ if it is right-handed.

\item
$6_1$ is $D(4)$, so that  $\Phi(6_1)=9+18t^0=27$.

\item
$7_4$ is deformed to $T'(6)$ as depicted in Fig.~\ref{74} by a smoothing
at a crossing where all colors are the same. Note that the colors
indicated in the figure exhaust all possibilities as $7_4$ has
cyclic Alexander module (by the result of Inoue \cite{Inoue}).
Hence we find $\Phi(7_4)=9 + 18t$.

\item
$7_7$ is deformed to $T'(6)$ as depicted in Fig.~\ref{77},
hence  $\Phi(7_7)=9 + 18t$.
The first deformation is a
smoothing
the second is isotopy.
Compare the left bottom diagram with the second left entry of
Fig.~\ref{antipara}. The left and right half of the diagram in
Fig.~\ref{77} are identified with the figure in Fig.~\ref{antipara},
and hence can be replaced by the top left entry of
Fig.~\ref{antipara}. The result is $T'(6)$ as depicted in bottom right
entry in Fig.~\ref{77}.

\end{itemize}

  \end{remark}

\begin{figure}
\begin{center}
\mbox{
\epsfxsize=3in
\epsfbox{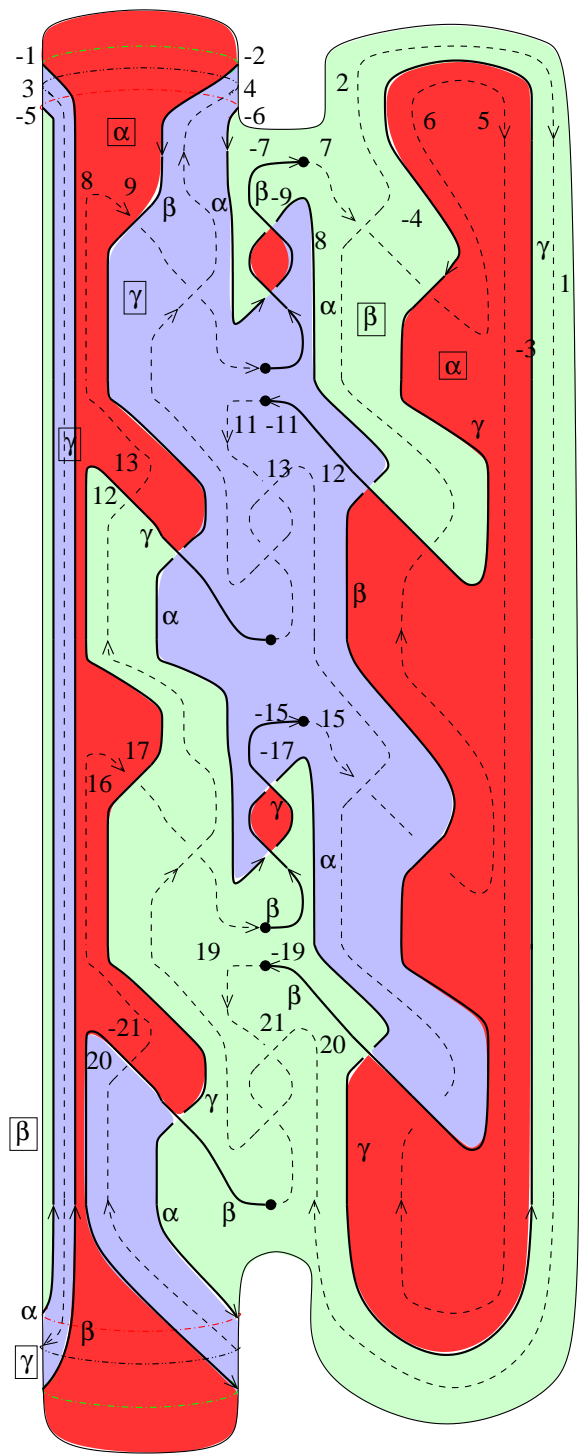}
}
\end{center}
\caption{The decker set of the $2$-twist spun trefoil}
\label{2twist}
\end{figure}

\begin{figure}
\begin{center}
\mbox{
\epsfxsize=2.7in
\epsfbox{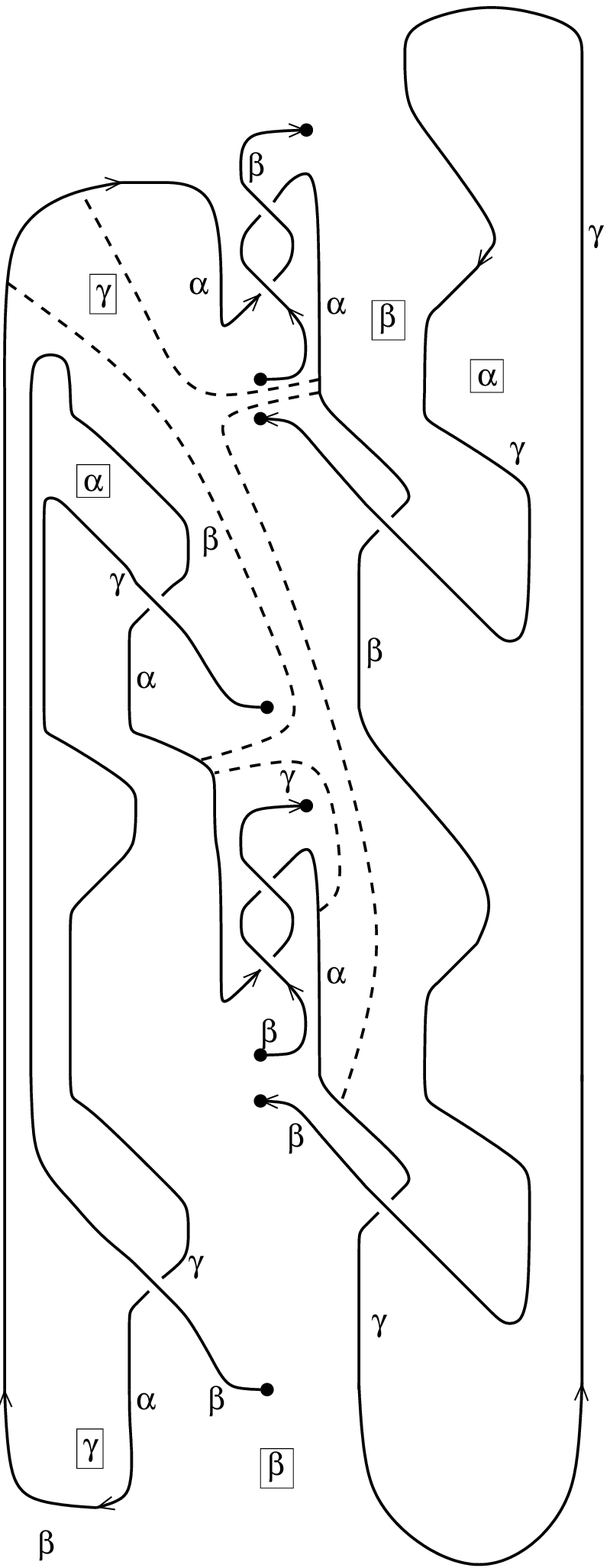}
}
\end{center}
\caption{The lower decker set of the $2$-twist spun trefoil}
\label{2twistlower}
\end{figure}

\bigskip

The computational technique discussed in this section can be applied to
knotted surfaces. We illustrate this with the $2$-twist spun trefoil.

\begin{example}{\rm Shadow coloring is closely related to
  coloring a knotted surface diagram and its lower decker set.
  The correspondence
was given in \cite{SSS2}, but we sketch the notion briefly.
 
Given an embedded surface in $4$-space, we chose a generic projection
  into $3$-space and label the double points of the projection above and
below to indicate their  relative distance  from the $3$-space
  into which they are projected. The pre-image of the double
  point set on the surface is called the {\it double decker set}.
  It is separated into  upper and lower pieces called
the {\it upper decker set} and the {\it lower decker set}.

A {\it quandle} coloring of a knotted surface diagram induces a shadow
  coloring of the lower decker set considered as an abstract arc diagram
  in the surface. The scheme for determining this coloring is as follows:
Along  a double point arc, an upper sheet with color $y$ locally separates
the lower sheet into two pieces (by the broken surface convention).
One component of the lower sheet is colored $x$ and the other
sheet is colored $x*y$
as depicted on the left of Fig.~\ref{triplepoint}.
The arc of lower decker points is colored  $y$
which represents the color of the corresponding over crossing sheet.
A branch point is an
  end-point of  a lower decker arc. At such a point,
  the colors on the arc and surrounding $2$-dimensional region coincide.
In the lowest of the three sheets that intersect at a triple point,
  two lower decker arcs intersect. The arc that corresponds to the upper
  most sheet is depicted as un-broken on the surface, and the middle arc
is broken into arcs locally. The colors on the arcs and regions
  match the shadow color condition depicted in
Fig.~\ref{3chain} (see also Fig.~\ref{triplepoint}).

In Figure ~\ref{2twist}, the double decker set for
the $2$-twist spun trefoil is depicted.
We give a brief description on how to obtain this diagram from a movie
of the $2$-twist spun trefoil. More details can be found in
\cite{CS:book}.
The solid lines correspond
to the lower decker points, and the dashed lines correspond to
the upper decker points.
The diagram has been colored by the $3$-element dihedral
  quandle $R_3=  \{0,1,2\}$  where $ \{\alpha,\beta,\gamma\} = \{0,1,2\}$.
The integer labels on double arcs are Gauss codes.
Even parity labels correspond to positive classical crossings,
and odd parity labels
correspond to negative crossings. Negative signs indicate undercrossings;
thus all solid lines have negative labels.
For example,
the first type II move in a movie of the knotted surface will
induce the birth of a negative  and a positive crossing that are labeled
$1$ and $2$, repsectively.
Immediately after the top saddle we have the Gauss code
$(-5,3,-1,-2,4,-6,2,-4,6,5,-3,1)$ which is the standard code
of the square knot.
The diagram for this decker set also indicates the height with respect
to the movie direction at which the critical points
occur.
For example the crossing in the center of the diagram that involves
arcs labeled $-7$, and $-9$ corresponds to the crossings between arcs
  labels $9$ and $2$ (on the left of the figure)
and between $7$ and $2$ on the right of the figure. All three of these
crossings correspond to a Reidemeister type III move among crossing
  points labeled $2$, $7$ and $9$.

In Fig.~\ref{2twistlower}, only the lower decker set and colors are depicted.

We now compute the state-sum invariant of this surface by the following 
technique.
Perform surgeries along the dotted arcs to
obtain four copies of the negative of the generator given in
Fig.~\ref{gentre}. Hence this color contributes $-1$ as the Kronecker
product with the pull-back of the $3$-cocycle $\xi$.
Hence we obtain $\Phi(K)=3+6 t^2$, as  we
computed
in \cite{CJKLS} by a different method.

} \end{example}

\section{Extensions of quandles by cocycles}

Let $X$ be a quandle and $A$ be an abelian group
written multiplicatively.
Let $\phi \in Z^2_{\rm Q} (X; A)$.
Consider $A$ as a trivial quandle ($a*b=a$ for any $a, b \in A$).
Let $E(X,A, \phi)$ be the quandle defined on the set $A \times X$ by
the operation $(a_1, x_1) * (a_2, x_2) = (a_1 \phi(x_1, x_2), x_1 * x_2)$.

\begin{lemma} \label{prodlemma}
The above defined operation $*$ on $A \times X$ indeed
defines a quandle $E(X,A, \phi) = (A \times X, *)$.
\end{lemma}
\begin{proof}
The idempotency is obvious.
For any $(a_2, x_2), (a, x) \in A \times X$,
let $x_1 \in X$ be a unique element $x_1 \in X$ such that
$x_1 * x_2 = x$. Then let $a_1 = a \phi(x_1, x_2)^{-1}$.
Then it follows that $(a_1, x_1) * (a_2, x_2) =(a, x) $,
and the uniqueness of  $(a_1, x_1)$ with this property is obvious.
The self-distributivity follows from
the $2$-cocycle condition by computation, as follows.
\begin{eqnarray*}
\lefteqn{
   [ (a_1, x_1) * (a_2, x_2)] * (a_3, x_3) } \\
  &=&
(a_1 \phi(x_1, x_2), x_1 * x_2) * (a_3, x_3) \\
&=& (a_1 \phi(x_1, x_2) \phi(x_1 * x_2,  x_3) ,  (x_1 * x_2)* x_3 ) ,
\end{eqnarray*}
and
\begin{eqnarray*}
\lefteqn{
  [ (a_1, x_1) * (a_3, x_3)] *[ (a_2, x_2) * (a_3, x_3)]} \\
  &=&
  (a_1 \phi(x_1, x_3), x_1 * x_3) *  (a_2 \phi(x_2, x_3), x_2 * x_3) \\
&=& (a_1 \phi(x_1, x_3) \phi( x_1 * x_3,   x_2 * x_3),
(x_1 * x_3) *(x_2 * x_3) ).
\end{eqnarray*}
We remark here that the computation above can be seen in knot diagrams
in Fig.~\ref{nu2cocy}. Go along the string that goes from
top left to bottom right in Reidemeister type III move
and read off the cocycles assigned at crossings. Then it picks up the
cocycles in the above two computations,
  for before and after the Reidemeister move, respectively.
\end{proof}

\begin{definition}
Two surjective homomorphisms of quandles $\pi_j: E_j \rightarrow X$,
$j=1,2$,
are called {\it equivalent } if there is
a quandle isomorphism $f: E_1 \rightarrow E_2$
such that $\pi_1=\pi_2 f$.
  \end{definition}

Note that there is a natural surjective homomorphism
$\pi: E(X,A, \phi)=A \times X  \rightarrow X$, which is the projection to
the second factor.

\begin{lemma}
If   $\phi_1$ and $\phi_2$ are cohomologous, i.e.,
$[\phi_1]=[\phi_2] \in  H^2_{\rm Q}(X;A)$,
then
$\pi_1: E(X,A, \phi_1) \rightarrow X$ and  $\pi_2: E(X,A, 
\phi_2)\rightarrow X$ are equivalent.
\end{lemma}
\begin{proof}
There is a $1$-cochain $\eta \in C^1_{\rm Q}(X;A)$ such that
$\phi_1 =\phi_2 \delta \eta$.
We show that
  $f: E(X,A, \phi_1)=A \times X \rightarrow A \times X = E(X,A, \phi_2)$
defined by $f(a,x)=(a \eta(x), x) $ gives rise to an equivalence.
First we compute
\begin{eqnarray*}
  f( (a_1, x_1) * (a_2, x_2) ) &  = & f( (a_1 \phi_1( x_1, x_2 ) ,  x_1 * 
x_2 ) ) \\
  & = &  (a_1 \phi_1( x_1, x_2 ) \eta( x_1 * x_2 ) ,  x_1 * x_2 ), \; \; 
\mbox{and} \\
  f( (a_1, x_1) ) * f(  (a_2, x_2) )  & = & (a_1 \eta(x_1) , x_1) * (a_2 
\eta (x_2) ,  x_2) \\
&=&  (a_1 \eta(x_1) \phi_2(x_1, x_2) ,   x_1 * x_2 )
\end{eqnarray*}
which are equal since $\phi_1 =\phi_2 \delta \eta$.
Hence $f$ defines a quandle homomorphism.
The map $f': A \times X \rightarrow A \times X$ defined by
$f'(a, x)=(a \eta(x)^{-1}, x)$ defines the inverse of $f$,
hence $f$ is an isomorphism. The map $f$ satisfies
$\pi_1=\pi_2 f$ by definition.
\end{proof}

\begin{lemma}
If
natural surjective homomorphisms 
(the projections to the second factor
$A \times X \rightarrow X$)      
$E(X,A, \phi_1)\rightarrow X$ and  $E(X,A, \phi_2)\rightarrow X$
  are equivalent, then  $\phi_1$ and $\phi_2$ are cohomologous:
$[\phi_1]=[\phi_2] \in  H^2_{\rm Q}(X;A)$.
\end{lemma}
\begin{proof}
Let $f:  E(X,A, \phi_1)=A \times X \rightarrow A \times X = E(X,A, \phi_2)$
be a quandle isomorphism with $\pi_1=\pi_2 f$.
Since $\pi_1(a,x)=x= \pi_2( f(a, x))$,
  there is an element $\eta(x) \in A$ such that
$f(a, x)= (a \eta(x), x)$,  for any $x \in X$.
This defines a function $\eta: X \rightarrow A $, $\eta \in C^1_{\rm Q}(X;A)$.
The condition that $f$ is a quandle homomorphism
implies that  $\phi_1 =\phi_2 \delta \eta$ by
the same computation as the preceding lemma.
Hence the result follows.
\end{proof}

The lemmas imply the following theorem.

\begin{theorem}
There is a bijection between the equivalence classes of
natural surjective homomorphisms 
  $E(X,A, \phi)\rightarrow X$ for a fixed $X$ and $A$,
   and
the set $H^2_{\rm Q}(X;A)$.
\end{theorem}

\begin{example}
The dihedral quandle $R_4$ is an abelian extension of $T_2$
by a non-trivial cocycle.
Let
$$\phi=\chi_{(0,1)} + \chi_{(1,0)} \in Z^2_{\rm Q}(T_2; {\bf Z}_2)$$
where  $\chi$ denotes the characteristic function
$$\chi_{x}(y) = \left\{ \begin{array}{ll} 1 & {\mbox{\rm if}} \  \ x=y
\\
                                                0   & {\mbox{\rm if}} \ \  x\ne
y  . \end{array} \right.$$
Let $X=X(T_2, \phi)$.
Then an isomorphism $f: R_4 \rightarrow X$ is defined by
$$ f(x) = ( \;  \lfloor x/2 \rfloor , \; \; \;  x \pmod{2} \;  ) , $$
where $\lfloor a  \rfloor $ denotes the largest integer not exceeding $a$.
It is checked directly that $f$ gives an isomorphism,
  but we can also see it as follows.
The quandle operation $a*b$ in $R_4$ is characterized by
the following property: if $a$ and $b$ have  the same parity,
then $a*b=a$. Otherwise, $a*b=\bar{a}$, where $\bar{a}$ has the same parity
as $a$ but is distinct from $a$.
For an element $(x, y) \in X={\bf Z}_2 \times T_2$, regard
$y$ as the parity. Then by the definition of $\phi$,
we see that $X$ has the quandle operation with the same property as above.
  \end{example}

\begin{example}
The quaternion group $Q_8=\{ \pm 1, \pm i , \pm j, \pm k \}$
under conjugation is a quandle. The set $\{ \pm 1 \}$ acts
  trivially on the rest as quandle operation,
and the rest $Q_6=\{ \pm i , \pm j, \pm k \}$ forms
a subquandle. The obvious projection $\pi: Q_6 \rightarrow T_3=\{ I, J, K\}$
defined by $\pi(\pm i)=I, \pi(\pm j)=J, \pi(\pm k)=K$
is a surjective quandle homomorphism.
One sees that $Q_6=E(X,A, \phi)$ where $X=T_3$, $A={\bf Z}_2$, and
$\phi \in Z^2_{\rm Q}(X;A)$ is defined by
$\phi = \prod_{a \neq b, \; a,b \in T_3 } \chi_{a, b}$.
  \end{example}

\begin{remark}
Let $\pi: Q\rightarrow X$ be a surjective quandle homomorphism,
such that the {\it equalizer} $E_{\pi}(x)=\{ \ z \in E \ | \ \pi(z)=x \ \} $
has the same (finite) cardinality.
Then it is an interesting problem to
determine when  $Q$
  is isomorphic to $E(X,A, \phi)$ for some $A$ and
  $\phi \in Z^2_{\rm Q} (X;A) $.
\end{remark}



\bigskip

Next we consider interpretations of $3$-cycles in extensions of quandles.
Let
$1 \rightarrow N   \stackrel{i}{\rightarrow} G  \stackrel{p}{\rightarrow} A 
\rightarrow 1$
be a short exact sequence of abelian groups.
Let $X$ be a quandle. For $\phi \in Z^2(X;A)$,
let $E(X,A, \phi)$ be as in the preceding section.
Let $s: A \rightarrow G$ be a set-theoretic
(not necessarily group homomorphism)
section,
i.e., $ps=\mbox{id}_A$ and $s(1_A)=1_G$.  

Consider the binary operation
$ (G \times X ) \times ( G \times X)
\rightarrow G \times X$ defined by
\begin{equation} \label{extdef}
  (g_1, x_1) * (g_2, x_2) = (g_1  s \phi (x_1, x_2) , x_1 * x_2 ).
\end{equation}
  We describe an obstruction to this being a quandle operation
by $3$-cocycles.

Since $\phi $ satisfies the $2$-cocycle condition,
$$p ( s \phi(x_1, x_2) s \phi (x_1 * x_2, x_3))
= p ( s \phi (x_1, x_3) s \phi (x_1 * x_3, x_2 * x_3 ) ) $$
  in $A$.
Hence there is a function
  $\theta: X \times X \times X \rightarrow N$ such that
\begin{equation} \label{3cocydef}
   s \phi (x_1, x_2)  s \phi (x_1 * x_2, x_3)
=i \theta (x_1, x_2, x_3)  s \phi (x_1, x_3) s \phi (x_1 * x_3, x_2 * x_3 ).
\end{equation}

\begin{lemma}
$\theta \in Z^3_{\rm Q}(X;N)$.
\end{lemma}
\begin{proof}
First, if $x_1=x_2$, or $x_2=x_3$, then the above defining relation 
for $\theta$ implies that $\theta (x_1, x_1, x_3)=1=\theta(x_1, x_2, x_2)$. 
For the $3$-cocycle condition, 
one computes 
\begin{eqnarray*}
\lefteqn{ \underline{ s \phi (x_1, x_2)  s \phi (x_1 * x_2, x_3)}
  s\phi ( (x_1 * x_2)*x_3, x_4) } \\
  &=& i \theta (x_1, x_2, x_3) [
  s \phi (x_1, x_3) \underline{ s \phi (x_1 * x_3, x_2 * x_3 ) } ]
  \underline{s\phi ( (x_1 * x_2)*x_3, x_4) } \\
  &=& [i\theta (x_1, x_2, x_3)i \theta (x_1 *x_3, x_2 * x_3, x_4) ]
\\ & &
[ \underline{ s \phi (x_1 * x_3,  x_4)} s\phi ((x_1 *x_3)*x_4,  (x_2 * x_3 
)*x_4 )]  \underline{ s \phi (x_1, x_3)}
\\ &=&
  [i \theta (x_1, x_2, x_3)i \theta (x_1 *x_3, x_2 * x_3, x_4)
i\theta(x_1, x_3, x_4)]\\ & &
  [s\phi (x_1,  x_4) s\phi (x_1*x_4,  x_3*x_4) ]
  s\phi ((x_1 *x_3)*x_4,  (x_2 * x_3 )*x_4 )
\end{eqnarray*}
and on the other hand,
\begin{eqnarray*}
\lefteqn{ s \phi (x_1, x_2)
\underline{  s \phi (x_1 * x_2, x_3)
  s\phi ( (x_1 * x_2)*x_3, x_4) } }\\
  &=& \underline{  s \phi (x_1, x_2) }
[ i\theta (x_1 * x_2, x_3, x_4) \underline{ s \phi (x_1 * x_2, x_4) }
  s\phi ( (x_1 * x_2)*x_4, x_3*x_4 ) ] \\
  &=&  i\theta (x_1 * x_2, x_3, x_4)
  \underline{   s\phi ( (x_1 * x_2)*x_4, x_3*x_4 )} \\ & &
[i\theta( x_1, x_2, x_4)  s \phi (x_1, x_4)
  \underline{   s \phi (x_1 * x_4 ,  x_2* x_4) } ] \\
  &=&  i\theta (x_1 * x_2, x_3, x_4) i\theta( x_1, x_2, x_4) s \phi (x_1, x_4)
[  i\theta(x_1 * x_4, x_2 * x_4, x_3 * x_4) \\ & &
  s \phi (x_1 * x_4 ,  x_3 * x_4)
  s \phi ( (x_1*x_4) * (x_3*x_4),  (x_2  * x_4)*( x_3 * x_4) ) ]
\end{eqnarray*}
so that we obtain the result.
The underlines in the equalities indicates where the relation \ref{3cocydef}
is going to be applied in the next step of the calculation.
\end{proof}

\begin{figure}
\begin{center}
\mbox{
\epsfxsize=3.5in
\epsfbox{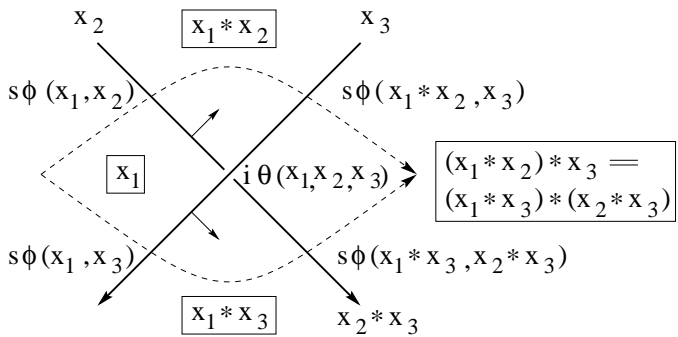}
}
\end{center}
\caption{Paths in shadow diagrams and $3$-cocycles  }
\label{shadowpath}
\end{figure}

\begin{figure}
\begin{center}
\mbox{
\epsfxsize=3.5in
\epsfbox{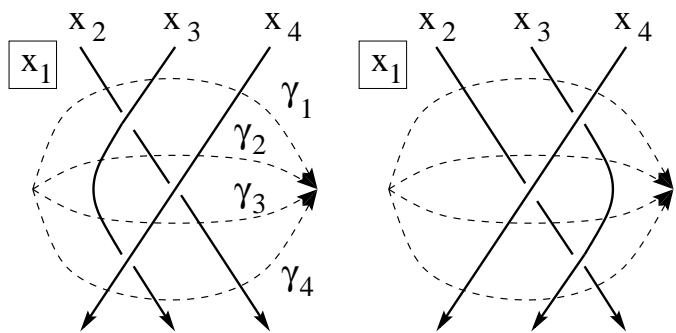}
}
\end{center}
\caption{The type III move with shadow colors and paths   }
\label{shadowIII}
\end{figure}

\begin{remark}
The above calculations are based on the knot diagrams as follows.
In Fig.~\ref{shadowpath}, two paths (denoted by dotted arcs)
  going from the left region
to the right region are depicted in a shadow colored diagram
near a crossing.
Let $\gamma_1$ be the top path, $\gamma_2$ be the bottom path.
Assign the cocycle $s\phi(x_1, x_2)$ when a dotted arc crosses
an over-arc with color $x_2$ from the region colored $x_1$ to
the region colored $x_1*x_2$ in the same direction as
the normal to the over-arc.
The product of these cocycles along $\gamma_1$ is
the left-hand side of Equality (\ref{3cocydef}), and
the product for $\gamma_2$ appears in the right-hand side.
The corresponding $3$-cocycle $i \theta(x_1, x_2, x_3)$
is assigned to the crossing.
We see this situation as follows. As we homotop
$\gamma_1$ to $\gamma_2$ through a crossing, the descrepancy of
the products is the $3$-cocycle assigned at the crossing.

In Fig.~\ref{shadowIII}, the shadow colors before/after
the Reidemeister type III move is specified. The paths denoted by
dotted arcs are also depicted, and named $\gamma_i$, $i=1,2,3,4$.
As we homotope the paths from $\gamma_1$ to $\gamma_4$,
we read off the $3$-cocycles assigned to the crossings.
The left and right of the Reidemeister move correspond to
the above two computations in the proof, respectively.

\end{remark}

Let $s': A \rightarrow G$ be another section, and  $\theta ' $
be a $3$-cocycle defined similarly for $s'$ by
\begin{equation} \label{3cocydef2}
   s' \phi (x_1, x_2)  s' \phi (x_1 * x_2, x_3)
=i \theta ' (x_1, x_2, x_3)  s' \phi (x_1, x_3) s' \phi (x_1 * x_3, x_2 * 
x_3 ).
\end{equation}

\begin{lemma}
The two $3$-cocycles $\theta$ and $\theta'$ are cohomologous,
$[\theta]=[\theta'] \in H^3_{\rm Q}(X;N)$.
\end{lemma}
\begin{proof}
Since $s(a)^{-1} s'(a) \in i (N) $
for any $a \in A$,
there is a function $\sigma: A \rightarrow N$ such that
$s'(a) =s(a) i \sigma(a)$ for any $a \in A$.
{}From Equality~(\ref{3cocydef2}) we obtain
\begin{eqnarray*}
\lefteqn{
   s \phi (x_1, x_2)  i \sigma \phi (x_1, x_2)
  s \phi (x_1 * x_2, x_3)  i \sigma\phi (x_1 * x_2, x_3) } \\
&=&
i \theta ' (x_1, x_2, x_3)  s \phi (x_1, x_3) i \sigma \phi (x_1, x_3)
  s \phi (x_1 * x_3, x_2 * x_3 ) i \sigma \phi (x_1 * x_3, x_2 * x_3 ) .
\end{eqnarray*}
Hence we have $\theta'= \theta \delta (\sigma \phi)$.
\end{proof}

\begin{lemma}
If $\theta$ is a coboundary, i.e.,
$[\theta ] =0 \in  H^3_{\rm Q}(X;N)$, then
$G \times X$ admits a quandle structure
such that $p \times \mbox{id}_X: G \times X \rightarrow A \times X$
is a quandle homomorphism.
\end{lemma}
\begin{proof}
By assumption there is $\xi \in C^2_{\rm Q}(X;N)$ such that
$\theta=\delta \xi$.
Define a binary operation on $G \times X$ by
$$
(g_1, x_1)*(g_2, x_2) = (g_1 \ s \phi(x_1, x_2) \xi (x_1, x_2)^{-1},
x_1 * x_2 ). $$
Then by Equality~(\ref{3cocydef}),
this defines a desired quandle operation.
\end{proof}

We summarize the above lemmas as

\begin{theorem} \label{obstthm}
The obstruction to extending  the quandle $E(X,A, \phi)=A \times X$ to
$G \times X$ lies in $H^3_{\rm Q}(X;N)$.
In particular, if $H^3_{\rm Q}(X;N)=0$, then $E(X,A, \phi)$
extends to $G \times X$ for any $\phi \in Z^2_{\rm Q}(X;A)$.
  \end{theorem}

\begin{sloppypar}
\begin{example}
For $X={\bf Z}_2[T, T^{-1} ]/(T^2+T+1)$, it is known
\cite{CJKLS} that $H^2_{\rm Q}(X;{\bf Z}_2)={\bf Z}_2$,
so that there is a non-trivial extension of $X$ to
$X({\bf Z}_2, \phi)= {\bf Z}_2 \times X$
defined by a generating $2$-cocycle $\phi$.
It is known that  $H^3_{\rm Q}(X;{\bf Z}_2)={\bf Z}_2^3$,
but we do not know which (if any) $3$-cocycles obstruct
extensions to $G \times X$ with
$1 \rightarrow {\bf Z}_2  \rightarrow G  \rightarrow {\bf Z}_2  \rightarrow 1.$
However, we obtain the following information from Theorem~\ref{obstthm}:
for any  $N \neq {\bf Z}_2$ with
$1 \rightarrow N  \rightarrow G  \rightarrow {\bf Z}_2  \rightarrow 1$,
the above non-trivial extension $E(X,A, \phi)$ further
extends to $G \times X$, as $H^3_{\rm Q}(X;N)=0$.
  \end{example}
\end{sloppypar}





\begin{thebibliography}{99}



\bibitem{Brieskorn} Brieskorn, E.,
{\it Automorphic sets and singularities,}
Contemporary math., 78 (1988), 45--115.


\bibitem{Brown}
Brown, K. S.,
{\it  Cohomology of groups.}
  Graduate Texts in Mathematics,
87. Springer-Verlag, New York-Berlin, 1982.

\bibitem{CJKLS}
  Carter, J.S.; Jelsovsky, D.; Kamada, S.; Langford, L.; Saito, M.,
{\it Quandle cohomology and state-sum invariants
of knotted curves and surfaces,}
preprint at
\begin{verbatim} http://xxx.lanl.gov/abs/math.GT/9903135 .\end{verbatim}

\bibitem{CJKS1}
  Carter, J.S.; Jelsovsky, D.; Kamada, S.; Saito, M.,
{\it Computations of quandle cocycle invariants of
  knotted curves and surfaces,}
preprint at
\begin{verbatim} http://xxx.lanl.gov/abs/math.GT/9906115 .\end{verbatim}

\bibitem{CJKS2}
  Carter, J.S.; Jelsovsky, D.; Kamada, S.; Saito, M.,
{\it Quandle homology groups, their betti numbers, and virtual knots,}
to appear in J. of Pure and Applied Algebra.


\bibitem{SSS2}  Carter, J.S.;  Kamada, S.; Saito, M.,
{\it Geometric interpretations of quandle
homology,}
to appear in Journal of Knot Theory and its Ramifications.


\bibitem{CS:book} Carter, J.S.; Saito, M.,
{\it Knotted surfaces and their diagrams,}
the American Mathematical Society,  1998.



\bibitem{FR}   Fenn, R.; Rourke,  C.,
\textit{Racks and links in codimension two,}
Journal of Knot Theory and Its Ramifications Vol. 1 No. 4 (1992), 343-406.



  \bibitem{FRS1}
Fenn, R.; Rourke, C.; Sanderson, B., {\it Trunks and classifying spaces,}
Appl. Categ. Structures 3 (1995), no. 4, 321--356.


\bibitem{FRS2} Fenn, R.; Rourke, C.; Sanderson, B.,
{\it James bundles and applications,} preprint found at
\begin{verbatim} http://www.maths.warwick.ac.uk/~bjs/ .\end{verbatim}

\bibitem{Flower} Flower, Jean, {\it Cyclic Bordism and Rack Spaces,} Ph.D. 
Dissertation, Warwick (1995).

\bibitem{FoxTrip} Fox, R.H.,
{\it A quick trip through knot theory,}
in Topology of $3$-Manifolds,
Ed. M.K. Fort Jr., Prentice-Hall (1962) 120--167.



\bibitem{GS}
Gerstenhaber, M.; Schack, S. D.,
{\it  Bialgebra cohomology,
deformations, and quantum groups.}
  Proc. Nat. Acad. Sci. U.S.A. 87 (1990), no. 1, 478--481.


\bibitem{Greene} Greene, M. T. {\it Some Results in Geometric Topology and 
Geometry,} Ph.D. Dissertation, Warwick (1997).

\bibitem{Inoue} Inoue, A.,
{\it Quandle homomorphisms of knot quandles to Alexander quandles,}
Preprint.


\bibitem{Joyce} Joyce, D.,
{\it A classifying invariant of knots, the knot quandle,}
J. Pure Appl. Alg., 23, 37--65.


\bibitem{K&P}  Kauffman, L. H., {\it Knots and Physics},
World Scientific, Series on knots and everything, vol. 1, 1991.

\bibitem{MS}
Markl, M.; Stasheff, J. D.
{\it  Deformation theory via deviations.}
J.
Algebra 170 (1994), no. 1, 122--155.


\bibitem{Matveev} Matveev, S.,
{\it Distributive groupoids in knot theory,} (Russian) Mat. Sb. (N.S.)
119(161) (1982), no. 1, 78--88, 160.

\bibitem{Rose} Roseman, D.,
{\it Reidemeister-type moves for surfaces in four dimensional space, }
in Banach Center Publications 42 (1998) Knot theory, 347--380.

\bibitem{RS} Rourke, C., and Sanderson, B.,
{\it There are two $2$-twist-spun trefoils,}
preprint at
\begin{verbatim} http://xxx.lanl.gov/abs/math.GT/0006062 .\end{verbatim}


\end{thebibliography}
\end{document}